\documentclass{amsart}
\usepackage{amssymb}
\newtheorem{theorem}{Theorem}[section]
\newtheorem{lemma}[theorem]{Lemma}
\newtheorem{corollary}[theorem]{Corollary}
\theoremstyle{definition}
\newtheorem{definition}[theorem]{Definition}

\numberwithin{equation}{section}

\def\IC{{\mathbb C}}
\def\IR{{\mathbb R}}
\def\cL{{\mathcal L}} 
\def\cH{{\mathcal H}} 
\def\bV{{\bf V}}
\def\diag{{\rm diag}\,} 
\def\conv{{\rm conv}\,}

\begin{document}

\title[Canonical forms, higher rank numerical ranges]
{Canonical forms, higher rank numerical ranges, totally isotropic subspaces,
and matrix equations}

\author{Chi-Kwong Li}
\address{Department of Mathematics, College of William \& Mary,
Williamsburg, VA 23185}
\email{ckli@math.wm.edu}
\thanks{Research of Li was partially supported by an NSF grant and 
a HK RGC grant. He is an honorary professor of the University of Hong 
Kong.}

\author{Nung-Sing Sze}
\address{Department of Mathematics, University of Connecticut,
Storrs, CT 06269}
\email{sze@math.uconn.edu}
\thanks{}



\dedicatory{}

\commby{}

\begin{abstract}
Results on matrix canonical forms are used to give a complete 
description of  the higher rank numerical range
of matrices arising from the study of quantum error correction.
It is shown that the set can be obtained as the 
intersection of closed half planes (of complex numbers). 
As a result, it is always a convex set in $\IC$. 
Moreover, the higher rank numerical range of a normal matrix
is a convex polygon determined by the eigenvalues. 
These two consequences confirm
the conjectures of Choi et al.\ on the subject.
In addition,  the results are used to derive a formula
for the optimal upper bound for the dimension of a totally isotropic 
subspace of a square matrix, 
and verify the solvability of certain matrix equations. 
\end{abstract}

\maketitle

\section{Introduction}

Let $M_n$ be the algebra of $n\times n$ complex matrices.
In \cite{Cet1}, the authors introduced the notion of the 
{\it rank-$k$ numerical range} of 
$A \in M_n$ defined and denoted by
$$\Lambda_k(A) = \{ \lambda \in \IC: PAP = \lambda P \hbox{\rm
~~ for some rank-$k$ orthogonal projection } P\}$$
in connection to the study of quantum error correction; see \cite{Cet2}.
Evidently, $\lambda \in \Lambda_k(A)$ if and only if there is a unitary matrix 
$U \in M_n$ such that $U^*AU$ has $\lambda I_k$ as the leading principal 
submatrix. 
When $k = 1$, this concept reduces to the classical numerical range,
which is well known to be convex by the Toeplitz-Hausdorff theorem;
for example, see \cite{L} for a simple proof. In \cite{Cet} the authors
conjectured that $\Lambda_k(A)$ is convex, and reduced the convexity
problem to the problem of showing that $0 \in \Lambda_k(T)$ for 
$$T = \begin{pmatrix}I_k & X \cr Y & -I_k \end{pmatrix}$$
for arbitrary $X, Y \in M_k.$
They further reduced this problem to the existence of a Hermitian
matrix $H$ satisfying  the matrix equation
\begin{equation} \label{eq1}
I_k + MH + HM^* - HPH = H
\end{equation}
for arbitrary $M \in M_k$ and positive definite $P \in M_k$.
In \cite{W}, the author observed that equation (\ref{eq1})
can be rewritten as the continuous Riccati equation 
\begin{equation} \label{eq2}
HPH-H(M^*-I_k/2)-(M-I_k/2)H - I_k = 0_k,
\end{equation}
and existing results on Riccati equation 
will ensure its solvability;
for example, see \cite[Theorem 4]{LR}. 
This establishes the convexity of $\Lambda_k(A)$.

Denote by $\lambda_k(H)$ the $k$th largest eigenvalue of
the Hermitian matrix $H \in M_n$.
We will use results on canonical forms of
complex square matrices to show that
$$\Lambda_k(A) = \bigcap_{\xi \in [0, 2\pi)}
\left\{\mu\in\IC: e^{i\xi}\mu +  e^{-i\xi}\bar \mu \le 
\lambda_k(e^{i\xi}A + e^{-i\xi}A^*)\right\}.$$
Thus, $\Lambda_k(A)$ is the intersection of closed half planes
on the complex plane, and therefore a convex set. Furthermore,
specializing our result to normal matrices
confirms the conjecture in \cite{Cet0} asserting that
$$\Lambda_k(A) = 
\bigcap_{1 \le j_1 < \cdots < j_{n-k+1} \le n} \conv\{ \lambda_{j_1},
\dots, \lambda_{j_{n-k+1}}\}$$
if $A \in M_n$ is a normal matrix with eigenvalues
$\lambda_1, \dots, \lambda_n$.
In addition, from our results one can derive a formula 
for the optimal upper bound for the 
dimension of a totally isotropic subspace
of a square  matrix.
As shown in \cite{Cet}, the convexity of the higher rank numerical range
is closely related to the study of solvability of matrix equations.
Following the idea in \cite{Cet}, we study 
the solvability of certain matrix equations including those 
of the  form (\ref{eq1}), (\ref{eq2}) and 
\begin{equation} \label{eq3}
I_k + RZ + Z^*R^* - Z^*Z = 0_k
\end{equation} 
for a given $k\times k$ matrix $R$.
In particular, it is shown that there is always a common solution 
$Z$ satisfying a pair of equations of the form (\ref{eq3}).
In other words, given two matrices $R, S \in M_k$,
the operator spheres
$$\{Z: |Z-R^*| = \sqrt{I_k+RR^*}\} \qquad 
\hbox{ and } \qquad  \{Z: |Z-S^*| = \sqrt{I_k+SS^*}\}$$
always have non-empty intersection,
here $|X|$ is the positive semidefinite square root of $X^*X$.

The following results on canonical forms of matrices 
will be used in our discussion; for 
example, see \cite{SS} and \cite{HS}.

\begin{itemize}
\item[I.] \it 
QR decomposition: \rm For every $A \in M_n$, there is a unitary matrix
$Q \in M_n$ and an upper triangular matrix $R\in M_n$ such that $A = QR$.
\item[II.] \it CS decomposition: \rm For every unitary $U \in M_{2k}$,  there
are unitary matrices $V = V_1 \oplus V_2$ and $W = W_1 \oplus W_2$
with $V_1, V_2, W_1, W_2 \in M_k$ such that
$$VUW = \begin{pmatrix}C & \sqrt{I_k - C^2}\cr 
\sqrt{I_k-C^2} & -C\cr\end{pmatrix},$$
where $C = \diag(c_1, \dots, c_k)$ with $1 \ge c_1 \ge \cdots \ge c_n \ge 0.$
\item[III.] \it 
${}^*$congruence canonical form: \rm
For every $A\in M_n$, there is an invertible 
$S\in M_n$ such that $S^*AS$ is a direct
sum of the following three types of matrices.
\begin{itemize}
\item[{\rm (1)}] 
$\Gamma_{2r}(\mu) = \begin{pmatrix}0_r & I_r\cr J_r(\mu) & 0_r \end{pmatrix} 
\in M_{2r}$ with $|\mu| > 1$ for $\mu \in \IC$,
where $J_r(\mu)$ is the $r\times r$ upper triangular Jordan block
with  eigenvalue $\mu$.

\item[{\rm (2)}]
$J_s(0)$, the $s\times s$ upper triangular Jordan block with eigenvalue zero. 

\item[{\rm (3)}] $e^{i\xi} \Delta_t$ with $\xi\in [0, 2\pi)$,
where $\Delta_t$ is the
$t\times t$ matrix whose $(p,q)$ entry
equals $1$ if $p+q = t+1$, equals the imaginary unit $i$ if $p+q=t+2$,
and equals 0 otherwise; in particular, $\Delta_1 = [1]$.
\end{itemize}

\end{itemize}

\section{Higher rank numerical range}
\setcounter{equation}{0}

\begin{definition}
For $A \in M_n$, let  
$\Omega_k(A)$ be the set of $\mu \in \IC$ such that
for each $\xi \in [0, 2\pi)$, the Hermitian matrix
$e^{i\xi}(A- \mu I_n) + e^{-i\xi}(A-\mu I_n)^*$ has at least 
$k$ nonnegative eigenvalues.
In particular, if $\lambda_k(H)$ denotes the $k$th largest 
eigenvalue of a Hermitian matrix $H \in M_n$,
then
$$\Omega_k(A) = \bigcap_{\xi \in [0, 2\pi)} 
\left\{\mu \in \IC: e^{i\xi}\mu+e^{-i\xi}\bar \mu 
\le \lambda_k(e^{i\xi} A+ e^{-i\xi}A^*) \right\}.$$
\end{definition}

\medskip
When $k = 1$, it is well known that the classical numerical range
$\Lambda_1(A)$ can be obtained by intersecting the closed half planes 
$$\{\mu \in \IC: e^{i\xi}\mu+e^{-i\xi}\bar \mu \le \lambda_1(e^{i\xi} A+ 
e^{-i\xi}A^*)\}, \qquad \xi \in [0, 2\pi).$$ 
We will show that $\Lambda_k(A) = \Omega_k(A)$,
which extends the classical result.
In particular, one can easily write a computer program to draw the boundary
$\partial \Omega_k(A)$ of $\Omega_k(A)$, and it is clear that 
for $A \in M_n$, the convex curve $\partial \Omega_k(A)$ lies inside 
the convex curve $\partial \Omega_{k-1}(A)$ if $k > 1$.

\medskip 
If $A$ is Hermitian, then we have the nested intervals
$$\Omega_1(A) \supseteq \Omega_2(A) \supseteq \Omega_3(A) \supseteq \cdots.$$
If $A$ is normal with eigenvalues 
$\lambda_1, \dots, \lambda_n$, then 
$$\Omega_k(A) = 
\bigcap_{1 \le j_1 < \cdots < j_{n-k+1} \le n} 
\conv\{ \lambda_{j_1}, \dots, \lambda_{j_{n-k+1}}\}$$ 
as described in \cite[Section 3]{Cet0}.
To see this, note that 
$$\mu \notin \bigcap_{1 \le j_1 < \cdots < j_{n-k+1} \le n} 
\conv\{ \lambda_{j_1}, \dots, \lambda_{j_{n-k+1}}\}$$ 
if and only if there is a line passing through 0 such that 
$n-k+1$ eigenvalues of $A-\mu I_n$ lie on one side
of the open half plane determined by the line; equivalently,
there is $\xi \in [0, 2\pi)$ such that $n-k+1$ eigenvalues of 
$e^{i\xi}(A- \mu I_n) + e^{-i\xi}(A- \mu I_n)^*$  are negative.

\medskip
Recall that we use $\lambda_k(H)$ to denote the $k$th largest eigenvalue of a
Hermitian matrix $H \in M_n$.
Our main theorem is the following.

\begin{theorem} \label{main} Let $A \in M_n$. Then 
$$\Lambda_k(A) = \Omega_k(A)
= \bigcap_{\xi \in [0, 2\pi)} 
\left\{\mu \in \IC: e^{i\xi}\mu+e^{-i\xi}\bar \mu 
\le \lambda_k(e^{i\xi} A+ e^{-i\xi}A^*) \right\}.$$
\end{theorem}

Since the intersection of half planes in $\IC$ is a convex set,
the following corollary is immediate.

\begin{corollary} \label{coro2.3}
Let $A \in M_n$. Then the rank-$k$ numerical 
range $\Lambda_k(A)$ is convex.
\end{corollary}

By the discussion on normal matrices before Theorem \ref{main},
we have the following corollary confirming the conjecture in \cite{Cet0}.

\begin{corollary} Let $A \in M_n$ be a normal matrix with eigenvalues 
$\lambda_1, \dots, \lambda_n$. Then 
$$\Lambda_k(A) =  
\bigcap_{1 \le j_1 < \cdots < j_{n-k+1} \le n} 
\conv\{ \lambda_{j_1}, \dots, \lambda_{j_{n-k+1}}\}.$$
\end{corollary}

To prove the theorem, we need the following lemma,
which  can be found in \cite{Cet}.
We give a short proof using the QR decomposition.

\begin{lemma} \label{lem0} Let $A \in M_n$ and $1 \le k \le n$.
Then $0 \in \Lambda_k(A)$ if and only if there is an invertible
$S \in M_n$ such that  $S^*AS$ has $0_k$ as the leading 
$k\times k$ principal submatrix.
\end{lemma}

\it Proof. \rm  The implication ``$\Rightarrow$'' is clear.
Conversely, suppose there is an invertible  $S \in M_n$ such that
$S^*AS$ has $0_k$ as the leading $k\times k$ principal submatrix.
By the QR decomposition, $S = UR$, 
where $U$ is unitary and $R$ is upper triangular.
Since $R^{-1}$ is also in upper triangular form, we see that 
$U^*AU = (R^{-1})^*(S^*AS)R^{-1}$ also has $0_k$ as 
the leading principal submatrix.
\qed

\medskip
We divide  the proof of Theorem \ref{main} into  three lemmas.
In particular, the construction in Lemmas \ref{lem2} and \ref{lem3} can
be done explicity using the results in \cite{FP,HS} and QR decomposition
(which involves only the Gram-Schmidt process). Thus, for every $\mu \in 
\Lambda_k(A)$, one can construct a unitary matrix $U$ such that $U^*AU$ with $\mu I_k$
as the leading principal submatrix.

\begin{lemma} \label{lem1}
Let $A \in M_n$.  Then $\Lambda_k(A) \subseteq \Omega_k(A).$
\end{lemma}

\it Proof. \rm Suppose $\mu \in \Lambda_k(A)$,
equivalently, $0 \in \Lambda_k(B)$
for $B = A - \mu I_n$. 
Then $e^{i\xi}B+e^{-i\xi}B^*$ 
is unitarily similar to 
a matrix with $0_k$ as the leading principal submatrix.
By the interlacing inequalities (for example, see  \cite{FP}),
$e^{i\xi}B+e^{-i\xi}B^*$ 
has at least $k$  nonnegative eigenvalues.
\qed

\begin{lemma} \label{lem2}
Let $A\in M_n$ be normal. Then $\Omega_k(A) \subseteq \Lambda_k(A)$.
\end{lemma}

\it Proof. \rm
Suppose $\mu \in \Omega_k(A)$. Let $B = A-\mu I_n$.
Then for each $\xi\in [0,2\pi)$,
the Hermitian matrix $e^{i\xi} B + e^{-i\xi}B^*$ 
has at least $k$ nonnegative eigenvalues.
We show that $0 \in \Lambda_k (B)$.

\medskip
We prove the result by induction on $k$.
If $k  = 1$, then the given condition ensures that 
$0$ lies in the convex hull of the eigenvalues of $B$. 
Suppose $V \in M_n$ is unitary such that
$V^*BV = \diag(b_1, \dots, b_n)$ and $p_1, \dots, p_n$
are nonnegative real numbers summing up to 1 such that
$\sum_{j=1}^n p_j b_j = 0$. 
Then $u = V(\sqrt{p_1}, \dots, \sqrt{p_n})^t$ is a unit
vector such that $v^*Bv = 0$. Choose a unitary matrix $U\in M_n$
with $u$ as the first column. Then
$U^*BU$ has zero as the $(1,1)$ entry.
So, the result holds for $k = 1$.
[One can also use the convexity of the classical numerical range to get the
conclusion. We include the argument so that the proof is 
independent of other convexity result.]

\medskip
Assume that $k > 1$ and the result is valid
for the rank-$m$ numerical range
of normal matrices whenever $m < k$. 
If $B$ has an eigenvalue equal to 0, then there is a unitary $V \in M_n$
such that $V^*BV = [0] \oplus B_1$ so that 
$e^{i\xi}B_1 + e^{-i\xi}B_1^*$ has at least $k-1$ 
nonnegative eigenvalues for any $\xi \in [0, 2\pi)$.
By induction assumption, there is a unitary $U \in M_{n-1}$
such that $U^*B_1U$ has $0_{k-1}$ as the leading principal submatrix.
Then $0_{k}$ will be a leading principal submatrix of
$([1]\oplus U)^*V^*BV([1]\oplus U)$. Thus, $0 \in \Lambda_k(B)$.

Now, assume that $B$ is invertible. Then $k \le n/2$. 
Suppose there is a pair of eigenvalues of $B$, say 
$\lambda_1$ and $\lambda_2$, satisfying 
$\lambda_1/|\lambda_1| = e^{i\theta}$
and $\lambda_2/|\lambda_2| = -e^{i\theta} = e^{i\theta + \pi}$ 
for some $\theta\in [0,2\pi)$. 
Then there is a unitary $V\in M_n$ 
such that $V^*BV = B_1 \oplus B_2$ with $B_1 = \diag(\lambda_1,\lambda_2)$.
Note that for each $\xi \in [0, 2\pi)$,
$e^{i\xi} B_1 + e^{-i\xi} B_1^*$ has at least $1$
nonnegative eigenvalue and
$e^{i\xi} B_2 + e^{-i\xi} B_2^*$ has at least $k-1$
nonnegative eigenvalues.
By the induction assumption, 
there are unitary $U_1 \in M_2$ and $U_2 \in M_{n-2}$
such that $U_1^*B_1U_1$ and $U_2^*B_2U_2$ have $0_1$ and $0_{k-1}$ 
as their leading principal submatrices, respectively.
Let $U = U_1 \oplus U_2$. 
Then $0_{k}$ will be a principal submatrix of
$U^* V^*BVU$ lying in rows and columns
$1,3,4 \dots, k+1$. Thus, $0 \in \Lambda_k(B)$.

\medskip
Continue to assume that $B$ is invertible; assume in addition that  
no pair of eigenvalues of $B$ 
have arguments $\theta$ and $\theta + \pi$.

\medskip\noindent
{\bf Claim} \it There is an invertible $S\in M_n$ 
such that $S^*BS$ has $0_{k}$ as the leading principal submatrix. \rm

\medskip
Once the claim is proved, we see that $0 \in \Lambda_k(B)$ 
by Lemma \ref{lem0}, and the induction proof will be complete.

To prove the claim, let
$\xi \in [0, 2\pi)$ be such that $e^{i\xi} B + e^{-i\xi}B^*$ 
has the smallest number of nonnegative eigenvalues, say, $k'$.
Then $k' \ge k$. We may assume that $k = k'$. 
Furthermore, we may assume that $\xi = 0$,
otherwise, replace $B$ by $e^{i\xi}B$.
Apply a ${}^*$-congruence to $B$ and assume that
$B = H+iG$ such that 
$H = I_k \oplus -I_{n-k}$
and $G = \diag(g_1, \dots, g_n)$ with
$g_1 \ge \cdots \ge g_k$ and $g_{k+1} \ge \cdots \ge g_n$.
Note that the given assumption on $B$ ensures that 

\begin{enumerate}
\item[\rm (i)] for every
straight line passing through the origin, there are at least $k$ eigenvalues
of $B$ lying in each of the closed half planes determined by the line,
and 
\item[\rm (ii)] there is no pair of eigenvalues of $B$ 
having arguments $\theta$ and $\theta + \pi$.
\end{enumerate}
\noindent
We claim that $-g_n > g_1$. Otherwise, $-g_n \le g_1$.
Since condition (ii) holds,  we see that $-g_n < g_1$. Moreover,
the line $\cL$ passing through
$0$ and the eigenvalue $1+ig_1$ of $B$ will divide the 
plane into two parts so that $k$ of the eigenvalues of $B$, namely,
$1+ig_1, \dots, 1+ig_k$ lies below $\cL$, and all other eigenvalues
lies in the open half plane above $\cL$. We may then rotate $\cL$
in the clockwise direction by a very small angle
so that at most $k-1$ of the eigenvalues of $B$,
namely, $1+ig_2, \dots, 1+ig_k$,
will lie on the closed half plane below the resulting line,
contradicting condition (i).

Similarly, we can argue that $-g_{n-1} > g_2$. Otherwise, $-g_{n-1} < g_2$,
and we can rotate the line passing through $0$ and $1+ig_2$ in the clockwise
direction by a very small angle so that at most $k-1$ eigenvalues
of $B$, namely, $1+ig_3, \dots, 1+ig_{k}$ and $-1+ig_n$,
will lie on the closed half plane below the resulting line, contradicting 
condition (i).

Repeating this argument, we see that
$$-g_n > g_1, \quad -g_{n-1} > g_2, \quad \cdots  \quad -g_{n-k+2} 
> g_{k-1},\quad  -g_{n-k+1} > g_k.$$
We can use a similar argument to show that
$$-g_{k+1} < g_k, \quad -g_{k+2} < g_{k-1},\quad \cdots \quad -g_{2k-1} 
< g_2, \quad -g_{2k} < g_1.$$
By \cite[Theorem 1]{FP}, there is a unitary 
$V \in M_{n-k}$ such that
$$V^*(\diag(g_{k+1}, \dots, g_n))V 
= \begin{pmatrix} -D & * \cr * & * \end{pmatrix}
\qquad \hbox{ with } D = \diag(g_1, g_2, \dots, g_k).$$
Thus, the leading $2k\times 2k$ submatrix of 
$(I_k \oplus V)^*B(I_k \oplus V)$ equals
$$(I_k+iD) \oplus (-I_k-iD).$$
Let 
$$W = \frac{1}{\sqrt 2}
\begin{pmatrix}I_k & I_k \cr I_k & - I_k\end{pmatrix}
\oplus I_{n-2k}.$$
Then the leading
$2k\times 2k$ submatrix of 
$W^*(I_k \oplus V)^*B(I_k \oplus V)W$ equals
$$\begin{pmatrix} 0_k & I+iD \cr I+iD & 0_k \end{pmatrix}.$$
So, the claim holds.
\qed

\begin{lemma} \label{lem3} For any matrix $A \in M_n$, 
we have $\Omega_k(A) \subseteq \Lambda_k(A)$.
\end{lemma}

\it Proof. \rm  Suppose $A \in M_n$ and $\mu \in \Omega_k(A)$.
Let $S \in M_n$ be such that $S^*(A - \mu I_n)S$ is
a direct sum of the following matrices as defined in Section {\rm 1 (III)}.
\begin{itemize}
\item[{\rm (a)}] $\Gamma_{2r_1}(\mu_1), \dots, \Gamma_{2r_u}(\mu_u)$.
\item[{\rm (b)}] $J_{s_1}(0), \dots, J_{s_v}(0)$, where $s_1, \dots, s_p$
are odd, and $s_{p+1}, \dots, s_v$ are even.
\item[{\rm (c)}] $e^{i\xi_1}\Delta_{t_1}, \dots, 
e^{i\xi_w}\Delta_{t_w}$, where $t_1, \dots, t_q$ are odd, and
$t_{q+1}, \dots, t_w$ are even.
\end{itemize}
Let $B = S^*(A-\mu I_n)S$.
For each $\xi \in [0, 2\pi)$, consider $e^{i\xi}B + e^{-i\xi}B^*$.
Each type (a) direct summand has the form $e^{i\xi}\Gamma_{2r_j}(\mu_j) + 
e^{-i\xi}\Gamma_{2r_j}(\mu_j)^*$, which will  contribute $r_j$ nonnegative (positive)
eigenvalues to $e^{i\xi}B+e^{-i\xi}B^*$.
Consequently, these summands will contribute
a total of $\sum_{j=1}^u r_j$ nonnegative
eigenvalues to $e^{i\xi}B+e^{-i\xi}B^*$.

\medskip
Each type (b) direct summand has the form 
$e^{i\xi}J_{s_j}(0) \oplus e^{-i\xi}J_{s_j}(0)^*$, which will contribute
$[(s_j+1)/2]$ nonnegative
eigenvalues to $e^{i\xi}B+e^{-i\xi}B^*$,
where $[x]$ denotes the integral part of the real number $x$.
Consequently, these summands will contribute a total of 
$\frac{1}{2} \left(\sum_{j=1}^v s_j + p\right)$ nonnegative
eigenvalues to $e^{i\xi}B+e^{-i\xi}B^*$.

\medskip
Each type (c) direct summand has the form
\begin{equation}
\label{tj}
e^{i(\xi+\xi_j)}\Delta_{t_j} + e^{-i(\xi+\xi_j)}\Delta_{t_j}^*
= \begin{pmatrix}
0 & & & & a_j \cr
& & & \cdot & b_j \cr
& & \cdot & \cdot & \cr
& \cdot & \cdot & & \cr
a_j & b_j & & & 0 \end{pmatrix}
\end{equation}
with $a_j = \cos (\xi+ \xi_j)$ and $b_j = -\sin (\xi + \xi_j)$.
Suppose $t_j$ is even. Since there is a $0_{t_j/2}$ leading principal
submatrix, the matrix has at least $t_j/2$ nonnegative eigenvalues.
If $\xi$ is chosen so that $a_j\ne 0$, then there will be
exactly $t_j/2$ nonnegative (positive) eigenvalues.
Thus, the matrix in (\ref{tj}) will contribute $t_j/2$ nonnegative
eigenvalues to $e^{i\xi}B+e^{-i\xi}B^*$.
Suppose $t_j$ is odd,  
then $e^{i(\xi+\xi_j)}\Delta_{t_j} + e^{-i(\xi+\xi_j)}\Delta_{t_j}^*$
is congruent to 
$[e^{i(\xi+\xi_j)}+e^{-i(\xi+\xi_j)}] \oplus D_j$ such that $D_j$ has
$(t_j-1)/2$ nonnegative eigenvalues.
Consequently, if $\xi$ is chosen so that $a_j \ne 0$ in (\ref{tj})
whenever $t_j$ is even,
then these summands will contribute a total of 
$\frac{1}{2}\left(\sum_{j=1}^w t_j - q\right) + \ell(\xi)$
nonnegative eigenvalues to $e^{i\xi}B+e^{-i\xi}B^*$, where
$\ell(\xi)$ is the number of nonnegative eigenvalues
of $e^{i\xi}N + e^{-i\xi}N^*$ with 
$$N = \diag(e^{i\xi_1}, \dots, e^{i\xi_q}).$$
Denote by $\nu(H)$ the number of nonnegative
eigenvalues of the Hermitian matrix $H$,
and 
$\ell' = \min\{\nu(e^{i\xi}N + e^{-i\xi}N^*): \xi \in [0, 2\pi)\}.$
Then there are infinitely many choices of $\xi$ which attain $\ell'$.
So, we may choose $\xi$ to attain $\ell'$
with the additional assumption that $a_j \ne 0$ in (\ref{tj})
whenever $t_j$ is even. 
Let 
\begin{eqnarray*}
k' 
&=& \sum_{j=1}^u r_j + \frac{1}{2}\left(\sum_{j=1}^v s_j + p\right) 
+ \frac{1}{2}\left(\sum_{j=1}^w t_j - q\right) + \ell' \\
&=& \min\{ \nu(e^{i\xi}B + e^{-i\xi}B^*): \xi \in [0, 2\pi)\}. \nonumber
\end{eqnarray*}
Then $k' \ge k$. Hence,
the conclusion that  $0\in \Lambda_{k}(A-\mu I_n)$
will follow
once we show that $0\in \Lambda_{k'}(A-\mu I_n)$.

By our assumption, $S^*(A - \mu I_n)S$ is a direct sum of the 
matrices listed in (a) -- (c). For each direct summand $\Gamma_{2r_j}$
in (a), it is clear that the leading principal submatrix is $0_{r_j}$.
Thus, these direct summands contain
a zero principal submatrix of dimension $\sum_{j=1}^u r_u$.

For each direct summand $J_{s_j}(0)$ in (b), the principal submatrix
lying in rows and columns indexed by odd numbers is a zero principal
submatrix. Thus, these direct summands contain a zero principal 
submatrix of dimension $(\sum_{j=1}^v s_j + p)/2$.

For each direct summand $e^{i\xi_j}\Delta_{t_j}$
in (c), if $t_j$ is even then the leading principal
submatrix is $0_{t_j/2}$; if $t_j$ is odd, then the 
leading principal submatrix is $0_{(t_j-1)/2}$.
Thus, these direct summands contain a zero principal submatrix
of dimension $(\sum_{j=1}^w t_j - q)/2$. 
Moreover, these direct summands 
are permutationally similar to a matrix $T$
with $0_t \oplus N$ as the $(t + q) \times (t+ q)$ 
leading principal submatrix,
where $t = (\sum_{j=1}^w t_j - q)/2$ and 
$N = \diag(e^{i\xi_1},\dots, e^{i\xi_q})$.
By Lemma \ref{lem2}, $0 \in \Lambda_{\ell'}(N)$. 
Thus, there is a unitary matrix $V\in M_q$
such that $V^*NV$ has  $0_{\ell'}$ as the principal submatrix.
Then $(I_t \oplus V \oplus I_{n-t-\ell'})^* T (I_t 
\oplus V \oplus I_{n-t-\ell'})$
has $0_{t+\ell'}$ as the leading principal submatrix.

Now combining all these zero principal submatrices yields
a zero principal submatrix of dimension
\begin{equation}\label{k'} 
k' = \sum_{j=1}^u r_j + \frac{1}{2}\left(\sum_{j=1}^v s_j + p\right) 
+ \frac{1}{2}\left(\sum_{j=1}^w t_j - q\right) + \ell'.
\end{equation}
The result follows. 
\qed

\section{Totally isotropic subspaces and matrix equations}
\setcounter{equation}{0}

Let $A \in M_n$. A subspace $\bV$ of $\IC^n$ is a {\it totally
isotropic subspace} of $A$ if
$x^*Ay = 0$ for any $x, y \in \bV$. Note that
$U\in M_n$ is unitary such that the first $k$ columns
of $U$ form a totally isotropic subspace of $A$ 
if and only if $U^*AU$ has $0_k$ as its leading principal submatrix.
One can also write $A = H+iG$ and dicuss the totally
isotropic subspace of the Hermitian matrix pair $(H,G)$, i.e.,
a subspace $\bV$ of $\IC^n$ such that $x^*Hy = 0 = x^*Gy$ for all 
$x, y \in \bV$.
It is clear that $A\in M_n$ has a totally isotropic subspace of
dimension $k$ if and only if $0 \in \Lambda_k(A)$.
By Theorem \ref{main}, we have the following.

\begin{theorem} Let $A \in M_n$. 
Denote by $\nu(H)$  the number of nonnegative eigenvalues of
the Hermitian matrix $H$.
Then 
\begin{eqnarray*}
&&\min\{ \nu(e^{i\xi}A + e^{-i\xi}A^*): \xi \in [0, 2\pi)\}\\
&=&\max\{\dim \bV : \bV \hbox{ is a totally isotropic subspace of } A \}.
\end{eqnarray*}
\end{theorem}

\medskip
Note that the quantity
$\min\{ \nu(e^{i\xi}A + e^{-i\xi}A^*): \xi \in [0, 2\pi)\}$
is equal to $k'$ in {\rm (\ref{k'})}, where the 
quantities $r_1, \dots, r_u, s_1, \dots, s_v$, etc. are 
determined by the canonical form 
of $A$ under ${}^*$-congruence as in the proof of Lemma 2.8
by putting $B = A-0I$.
By the result in \cite{HS}, one
can obtain the canonical form $S^*AS$
by a finite algorithm using exact arithmetic.

The authors of \cite{Cet} showed that the study of the
convexity of the higher rank numerical range can be reduced to 
verifying the following lemma, which follows readily from
Corollary \ref{coro2.3}.

\begin{lemma} \label{thm1} Let 
$A = \begin{pmatrix}I_k & X \cr Y & -I_k\end{pmatrix} \in M_{2k}$,
where $X, Y \in M_k$. Then there is a unitary $U \in M_{2k}$
such that $U^*AU$ has $0_k$ as the leading  
$k\times k$ principal submatrix.
\end{lemma}

\it Proof. \rm Since $1, -1 \in \Lambda_k(A)$, we see that
$0 \in \Lambda_k(A)$ by Corollary \ref{coro2.3}. \qed

\medskip
In \cite{Cet}, it was shown that
the existence of $U$ in Lemma \ref{thm1} is equivalent 
to the solvability of some matrix equations;
see \cite[Theorem 2.12]{Cet}.
In the next theorem, we will use Lemma \ref{thm1} and the
CS decomposition of matrices to prove the 
solvability of a number of matrix equations and system of matrix 
equations. The equations in (a), (d), (f) have been considered 
in \cite{Cet}. We give slightly different proofs of them.

We consider also other matrix equations. In particular,
assertion (c) of the theorem can be restated as
$$\{Z: |Z-R^*| = \sqrt{I_k+RR^*}\} \cap \{Z: |Z-S^*| = \sqrt{I_k+SS^*}\}
\ne \emptyset.$$

One can use the results in \cite{FP, HS} and QR decomposition to
construct the unitary matrix $U$ in Lemma \ref{thm1}. As a result,
one can give an explicit construction of the
solutions of the matrix equations (a) -- (d)
following our proof. 

It is easy to check that  
solvability of the equations in the theorem is equivalent to the 
existence of a unitary $U \in M_{2k}$ satisfying the conclusion of 
Lemma \ref{thm1}.

As suggested by Professor T. Ando, it is interesting and inspiring to
consider the scalar case of the statements and the proofs of 
the equations in the theorem.

\begin{theorem} \label{thm3}
Let $R, S, P, C \in M_k$ such that $P$ is positive
definite, $C$ is a strict contraction, and $\gamma \in \IR$.

\begin{itemize}
\item[{\rm (a)}] There is a $Z \in M_k$ such that
$$I_k + RZ + Z^*S^* - Z^*Z = 0_k.$$
\item[{\rm (b)}]
There is a $Z \in M_k$ such that 
$$I_k+RZ+Z^*R^*-Z^*Z = 0_k \quad \hbox{ and } \quad
  SZ+Z^*S^* = 0_k.$$
\item[{\rm (c)}]
There is a $Z \in M_k$ such that 
$$I_k+RZ+Z^*R^*-Z^*Z = 0_k \quad \hbox{ and } \quad
  I_k+SZ+Z^*S^*-Z^*Z = 0_k.$$
\item[{\rm (d)}] There is a Hermitian $H \in M_k$ such that
$$I_k + RH + HR^* - HPH = \gamma H.$$
\item[{\rm (e)}] There is a unitary $U \in M_k$ such that
$$SR^* + RS^* = SU\sqrt{I_k+RR^*} + \sqrt{I_k+RR^*}U^*S^*.$$
\item[{\rm (f)}] There is a unitary $U\in M_k$ 
and a Hermitian $H\in M_k$ such that
$$C = U + PH.$$ 
\end{itemize}
\end{theorem}

\it Proof. \rm
Consider the equation in (a).
By Lemma \ref{thm1},
there is a unitary $U \in M_{2k}$ such that 
$$U^*BU = \begin{pmatrix}0_k & * \cr * & * \end{pmatrix}
\quad \hbox{ with } \quad 
B = \begin{pmatrix} I_k & R \cr S^* & -I_k \end{pmatrix}.$$
By the CS decomposition,
there are unitary matrices
$V = V_1\oplus V_2$, $W = W_1\oplus W_2 \in M_{2k}$ with
$V_1, V_2, W_1, W_2 \in M_k$ such that 
$$(V_1\oplus V_2)U(W_1\oplus W_2)=  
\begin{pmatrix}C & \sqrt{I_k-C^2} \cr \sqrt{I_k-C^2}& -C \end{pmatrix}$$
where $C = \diag(c_1, \dots, c_k)$ with $1 \ge c_1 \ge \cdots \ge c_k \ge 0$.
Then $W(U^*BU)W^*$ also has $0_k$ as the leading 
$k\times k$ principal submatrix. Equivalently,
\begin{eqnarray*}
0_k&=& \begin{pmatrix} C & \sqrt{I_k-C^2} \end{pmatrix}
\, V^*BV\,
\begin{pmatrix} C \cr \sqrt{I_k-C^2} \end{pmatrix} \\
&=& C^2 + CV_1^*RV_2\sqrt{I_k-C^2} + \sqrt{I_k-C^2}V_2^*S^*V_1C - (I_k-C^2).
\end{eqnarray*}
Evidently, $c_k > 0$. Otherwise, the $(k,k)$ entry of the above
matrix is $-1$. Thus, we can multiply the above equation by 
$V_1C^{-1}$ on the left and $C^{-1}V_1^*$ on the right to get
$$0_k = I_k + RZ + Z^*S^* - Z^*Z \qquad
\hbox{ with } \quad Z = V_2\sqrt{I_k-C^2}(C^{-1}V_1^*).$$

To prove (b), let $\tilde R = R+S$ and $\tilde S = R-S$.
By (a), there is $Z \in M_k$ such that 
$$I_k+\tilde R Z+Z^*\tilde S^*-Z^*Z = 0_k.$$
Taking the Hermitian part and skew-Hermitian part of the 
above equation, we get the two equations in (b).

To prove (c), let $\tilde S = S-R$. By (b) there is $Z \in M_k$ such that
$$I_k+RZ+Z^*R^*-Z^*Z = 0_k 
\qquad \hbox{and}\qquad
\tilde SZ+Z^*\tilde S^*= 0_k.$$
Adding the two equations, we get
$I_k+SZ+Z^*S^*-Z^*Z = 0_k.$

To prove (d), we may assume that $\gamma = 0$. Otherwise,
replace $R$ by $R - \gamma I_k/2$.
Let $S = iP^{-1/2}$ and $\tilde R = RP^{-1/2}$. 
By (b) there is $Z \in M_k$ so that
$$I_k+\tilde RZ+Z^*\tilde R^*-Z^*Z = 0_k \qquad 
\hbox{ and } \qquad SZ+Z^*S^*= 0_k.$$
The second equation implies that  $Z = P^{1/2}H$ for some 
Hermitian $H$. Putting $Z = P^{1/2}H$ to the first equation, we have
$I_k + RH + HR^* - HPH=0_k$ as asserted. 

To prove (e), note that the first equation in (b) can be 
written as  $(Z^*-R)(Z-R^*) = I_k + RR^*.$ Thus, its solution 
has the form $Z = R^* - U \sqrt{I_k+RR^*}$ for some unitary $U \in M_k$.
Substituting this into the second equation in (b),
we get the desired conclusion. 

Finally, to prove (f), let $S = iTP^{-1}$ and $R = T C^*$
with $T = (I_k-C^*C)^{-1/2}$. Then $\sqrt{I_k+RR^*} = T$.
By (e), there is a unitary $U\in M_k$ such that
$$(SR^* - SU\sqrt{I_k+RR^*}) + (RS^*  - \sqrt{I_k+RR^*}U^*S^*) = 
0_k.$$
Hence, $iTP^{-1}(C-U)T = SR^* - SU\sqrt{I_k+RR^*} = iK$
for some Hermitian $K$. 
Take $H = T^{-1} KT^{-1}$.  The result follows.
\qed

\section{Infinite Dimensional Operators and Related results}

One can easily extend the definition of $\Lambda_k(A)$
to a bounded linear
operator $A$ acting on infinite dimensional Hilbert spaces $\cH$;
for example, see \cite{W}.  
Results on $\Lambda_k(A)$ for infinite dimensional operators
have been obtained in \cite{LPS2} 
including Theorem \ref{4.1} below.
For a self-adjoint operator $H$, we let
$$\lambda_k(H) = \sup \{\lambda_k(X^*HX): X: \IC^k \rightarrow \cH,
\ X^*X = I_k \}.$$

\begin{theorem} \label{4.1}
Let $A$ be a bounded linear operator acting on an infinite
dimensional Hilbert space $\cH$.
Then $\Lambda_k(A)$ is convex and its closure equals
$$\bigcap_{\xi \in [0, 2\pi)} \left\{\mu \in \IC: 
e^{i\xi}\mu + e^{-i\xi}\bar\mu \le \lambda_k(e^{i\xi}A + 
e^{-i\xi}A^*)\right\}.$$
\end{theorem}

An open question in \cite{Cet0} concerns the 
lower bound of $\dim \cH$ which ensures that
$\Lambda_k(A)$ is non-empty for every bounded linear
operator $A$ acting on $\cH$.  
The following result was proved in \cite{LPS1}  
that answers the above question.

\begin{theorem} Let $\cH$ be a Hilbert space, and let
$k$ be a positive integer. Then $\Lambda_k(A) \ne \emptyset$
for every bounded linear operator $A$ acting on $\cH$
if and only if $\dim \cH > 3k-3$.
\end{theorem}

\bigskip\noindent
\bf Acknowledgment
\rm

\medskip
We would like to thank the authors of \cite{Cet} and \cite{W} for
sending us their preprints. We also thank Professors T. Ando,
V. Bolotnikov, M.D. Choi, R.A. Horn, Y.T. Poon and L. Rodman 
for some helpful correspondence and discussion.

\bibliographystyle{amsplain}

\end{document}